\documentclass[11pt]{article}

\usepackage{amsmath}
\usepackage{amssymb}
\usepackage{amscd}
\usepackage{latexsym}
\usepackage{theorem}
\usepackage{epsfig}
\usepackage{psfrag}
\usepackage{amsfonts}
\usepackage{euscript}

\newtheorem{theorem}{Theorem}[section]
\newtheorem{lemma}[theorem]{Lemma}
\newtheorem{proposition}[theorem]{Proposition}

{\theorembodyfont{\rmfamily}
\theoremstyle{plain}
\newtheorem{definition}[theorem]{Definition}

\newtheorem{remark}[theorem]{Remark}
\newtheorem{note}[theorem]{Remark}
\newtheorem{assumption}{Assumption}
}

\numberwithin{equation}{section}
\numberwithin{figure}{section}

\input xy
 \xyoption{all}

\newcommand{\PP}{{\mathbb P}}

\renewcommand{\Im}{\operatorname{Im}}

\newcommand{\Ker}{\operatorname{Ker}}
\newcommand{\Hom}{\operatorname{Hom}}

\newcommand{\C}{{\mathbb{C}}}

\newcommand{\Z}{{\mathbb{Z}}}

\newcommand{\R}{{\mathbb{R}}}
\newcommand{\N}{{\mathbb{N}}}

\newcommand{\G}{{\mathcal{G}}}
\renewcommand{\P}{{\mathcal{P}}}
\newcommand{\B}{{\mathcal{B}}}

\newcommand{\into}{{\hookrightarrow}}

\newcommand{\qed}{\hfill \mbox{$\Box$}\medskip\newline}
\newenvironment{proof}{\noindent {\bf Proof:}}{\qed \par}

\newcommand{\algt}{\mathfrak{t}}

\begin{document}

%\begin{spacing}{1.1}

%%%%%%%%%%%%%%%%%%%%
% Title
%%%%%%%%%%%%%%%%%%%%

\noindent {\LARGE \bf Computation of generalized equivariant
cohomologies\\ \medskip of Kac-Moody flag varieties }
\bigskip\\

%%%%%%%%%%%%%%%%%%%%
% Authors
%%%%%%%%%%%%%%%%%%%%

\noindent {\bf Megumi Harada }\footnote{{\tt megumi@math.toronto.edu}} \\
Department of Mathematics, University of Toronto,
Toronto, Ontario M5S 3G3 Canada\smallskip \\
{\bf Andr\'e Henriques }\footnote{{\tt andrhenr@math.mit.edu}}\\
Department of Mathematics, Massachusetts Institute of Technology,
Cambridge, MA 02139-4307 \smallskip\\
{\bf Tara S. Holm }\footnote{{\tt tsh@math.berkeley.edu} \\
\newline \mbox{~~~~}{\it MSC 2000 Subject Classification}:
Primary: 55N91  \hspace{0.1in} Secondary: 22E65, 53D20.
\newline \mbox{~~~~}{\it Keywords}:
equivariant cohomology, equivariant $K$-theory, equivariant complex cobordism, flag varieties, Kac-Moody groups, stratified spaces. \newline
\newline An earlier version of this paper, entitled {\em $T$-equivariant cohomology of cell 
complexes and the case of infinite Grassmannians}, is still available at 
{\bf math.DG/0402079.}}\\
Department of Mathematics, University of California, Berkeley, CA
94720-3890

%\vfill\pagebreak

%%%%%%%%%%%%%%%%%%%%
% Disclaimer
%%%%%%%%%%%%%%%%%%%%
%\begin{center}
%\framebox{
%{\Large\bf DRAFT (\today): DO NOT DISTRIBUTE.}}
%\end{center}

%%%%%%%%%%%%%%%%%%%%%
%  Abstract
%%%%%%%%%%%%%%%%%%%%%
{\small
\begin{quote}
\noindent {\em Abstract.} In 1998, Goresky, Kottwitz, and
MacPherson showed that for certain projective varieties $X$ equipped
with an algebraic action of a complex torus $T$, the equivariant
cohomology ring $H^*_T(X)$ can be described by combinatorial data
obtained from its orbit decomposition. In this paper, we generalize
their theorem in three different ways.  First, our group $G$ need not
be a torus. Second, our space $X$ is an equivariant stratified space,
along with some additional hypotheses on the attaching maps.  Third,
and most important, we allow for generalized equivariant cohomology
theories $E_G^*$ instead of $H_T^*$. For these spaces, we give a
combinatorial description of $E_G^*(X)$ as a subring of $\prod
E_G^*(F_i)$, where the $F_i$ are certain invariant subspaces of
$X$. Our main examples are the flag varieties $\G/\P$ of Kac-Moody
groups $\G$, with the action of the torus of $\G$. In this context,
the $F_i$ are the $T$-fixed points and $E_G^*$ is a $T$-equivariant
complex oriented cohomology theory, such as $H_T^*$, $K_T^*$ or
$MU_T^*$. We detail several explicit examples.
\end{quote}
}
\bigskip

\section{Introduction and Background}\label{intro}

The goal of this paper is to give a combinatorial description of
certain generalized equivariant cohomologies of stratified spaces.
The important examples to which our main theorems apply include
$T$-equivariant cohomology, $K$-theory, and complex cobordism
of Kac-Moody flag varieties.  Although the examples that motivate
us come from the theory of algebraic groups, our proofs rely
heavily on techniques from algebraic topology.  Indeed, we state
the results of Sections~\ref{se:injectivity} through
\ref{se:generators} in the following context.

Let $G$ be a topological group and $E_G^*$ a $G$-equivariant
cohomology theory (see \cite[Chapter~XIII]{May} for a definition)
with a commutative cup product. Let $X$ be a stratified $G$-space
such that successive quotients $X_i/X_{i-1}$ are homeomorphic to
Thom spaces $Th(V_i)$ of $E$-orientable $G$-vector bundles $V_i\to
F_i$. In this setting, and with the assumption that the Euler
classes $e(V_i)$ are not zero divisors, we show that the
restriction map
$$
\imath^*:E_G^*(X) \to \prod_i E_G^*(F_i)
$$
is injective.  Moreover, when $X$ and the $G$-action satisfy
additional technical assumptions, we identify the image of
$\imath^*$ as a subring of $\prod_i E_G^*(F_i)$ defined by
explicit compatibility conditions involving  divisibility by
certain Euler classes.  We also construct free $E_G^*$-module
generators of $E_G^*(X)$.

\begin{figure}[h]
\psfrag{x}{$x_1$} \psfrag{y}{$x_2$} \psfrag{x+y}{$x_1+x_2$}
\psfrag{2y}{$2x_2$} \psfrag{x+2y}{$x_1+2x_2$}
\begin{center}
\epsfig{figure=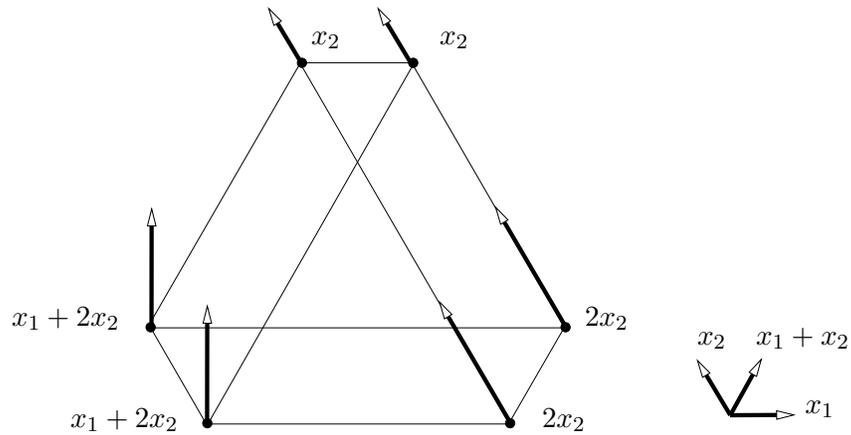,width=5in}
\end{center}
\begin{center}
\parbox{5in}{
\caption[coadjoint orbit]{This shows the graph $\Gamma$ for a flag
variety $SL(3,\C)/\B$.  The weight $\alpha_{(v,w)}$ is exactly the
direction of the edge $(v,w)$, as explained in
Section~\ref{sec:eqvtCohAffGr}. There is a linear polynomial
attached to each vertex, also depicted as a vector. 
The polynomials satisfy the compatibility
conditions, so this does represent an equivariant cohomology class
in $H_T^2(SL(3,\C)/\B)$.}\label{fig:SU3} }
\end{center}
\end{figure}

Our theorems generalize known results in algebraic and symplectic
geometry. When $X$ is a projective variety, $G$ a complex torus,
and $E_G^*$ ordinary equivariant cohomology, then we recover a
theorem of Goresky, Kottwitz and MacPherson \cite{GKM} that
computes $H_T^*(X;\C)$. They assume that $X$ has finitely many
$0$- and $1$-dimensional $T$-orbits, and then consider the graph
$\Gamma$ whose vertices are the fixed points $X^T$ and edges are
the one-dimensional orbits. An edge $(v,w)$ in $\Gamma$ is
decorated with the weight $\alpha_{(v,w)}$ of the $T$-action on
the corresponding orbit. They provide a combinatorial description
of $H_T^*(X)$ as a subring of $H_T^*(X^T)$ in terms of this graph.
Each edge of $\Gamma$ gives a condition as follows. Let $x(v)$
denote the restriction of a class $x\in H_T^*(X)$ to $v\in X^T$.
Then the condition reads
\begin{equation}\label{eq:intro}
\alpha_{(v,w)}\ \big|\ x(v) - x(w).
\end{equation}
We illustrate an example in Figure~\ref{fig:SU3}.

This article is organized as follows. In
Section~\ref{se:injectivity} we prove the injectivity of the map
$$
\imath^*:E_G^*(X) \to \prod_i E_G^*(F_i).
$$
Next, in Section~\ref{se:GKM}, we identify the image of
$\imath^*$, giving combinatorial conditions similar to those in
\eqref{eq:intro}. In Section~\ref{se:generators}, we give a
description of module generators for $E_G^*(X)$.  Finally, in
Sections~\ref{sec:eqvtCohAffGr} and \ref{sec:examples}, we return
to our motivating examples, which are homogeneous spaces $\G/\P$
for Kac-Moody groups $\G$, equipped with the action of a torus
$T$. For these spaces, our theory applies when $E_T^*$ is any
complex oriented $T$-equivariant cohomology theory. We make
explicit computations for three examples: a homogeneous space of
$G_2$, the based loop space $\Omega SU(2)$, and a homogeneous
space of \(\widehat{LSL(3,\C)}^{\Z/2\Z} \rtimes \C^*\).

\bigskip
\noindent{\bf \em Acknowledgments.} The first version of this
paper concerned ordinary $T$-equivariant cohomology and cell
complexes with even dimensional cells. We thank the referee for
pointing out that these results extend to arbitrary generalized
cohomology theories and more general stratified spaces. (S)he also
helped streamline our original proofs in
Section~\ref{se:injectivity}, and in particular offered a proof of
Theorem~\ref{thm:injectivity}.

We offer many thanks to the following people: to Allen Knutson for
originally suggesting the problem of a GKM theory for the
homogeneous spaces of loop groups and for teaching the first and
second authors how to draw GKM pictures; to Dylan Thurston for
providing a quick proof of Lemma~\ref{lemma:qqq};  to Dev Sinha
for helpful conversations about $MU_T^*$. The second author thanks
the Universities of Geneva and Lausanne for hosting and funding
him during part of this work. The third author was supported in
part by a National Science Foundation Postdoctoral Fellowship. All
authors are grateful for the hospitality of the Erwin
Schr\"odinger Institute in Vienna, where some of this research was
conducted.

\section{The injectivity theorem for stratified spaces}\label{se:injectivity}

Let $G$ be a topological group and $E^*_G$ a $G$-equivariant
cohomology theory with commutative cup product. We consider
stratified $G$-spaces
\begin{equation}\label{stratification}
X = \bigcup_{i \ge 1} X_i, \quad X_1 \subseteq X_2 \subseteq X_3 \ldots
\end{equation}
where the successive quotients \(X_i/X_{i-1}\) are homeomorphic to
the Thom spaces $Th(V_i)$ of some $G$-vector bundles \(V_i\to
F_i\). Moreover, we require that the above vector bundles be
$E$-orientable (see \cite[p.\ 177]{May}). In other words, X is
built by successively attaching disc bundles $D(V_i)$ via
equivariant attaching maps $\varphi_i:S(V_i)\to X_{i-1}$. This
should be compared to the way one builds CW complexes by
successively attaching discs.

We recall that an $E$-orientation, or Thom class, of a $G$-vector
bundle $V\to F$ is an element $u\in E^*_G(Th(V))$. For each closed
subgroup $H<G$ and point $x\in F^H$, the restriction of $u$ to
$V|_{G\cdot x}$ is a generator of the free $E^*_H$-module
$E^*_G(Th(V|_{G\cdot x}))\simeq E^*_H(D(V_x),S(V_x))$. The Euler
class $e(V)$ is the restriction of the Thom class $u$ to the base
$F$ via the zero section map.

\begin{note}\label{re:posets}
As with CW complexes, the stratification is often more naturally
indexed by a poset $I$ rather than $\mathbb N$. In that case, one
should replace the expression $X_i/X_{i-1}$ by
$X_i/\bigcup_{j<i}X_j$. The poset $I$ is required to satisfy the
condition that $\{j\in I:j<i\}$ is finite for all $i\in I$, which
makes the inductive proofs work. In the proofs, we ignore this
fact and pretend that $I=\N$. The only thing that we need is that
for each $i\in I$, the subspace $X_i$ is obtained by a finite
sequence of gluings, and that
$X=\displaystyle\underrightarrow{\lim}X_i$.
\end{note}

\begin{note}\label{wecare}
In the examples in Sections~\ref{sec:eqvtCohAffGr} and
\ref{sec:examples}, the group $G=T$ is a finite dimensional torus,
the $T$-spaces $F_i$ are single points and the $V_i$ are complex
$T$-representations. The stratification \eqref{stratification}
expresses $X$ as a cell complex with even dimensional cells.
\end{note}

The main theorem of this section establishes the injectivity of
the restriction map $E_G^*(X) \to E_G^*(\coprod F_i) \cong \prod
E_G^*(F_i)$ when the Euler classes are not zero divisors.

\begin{theorem}\label{thm:injectivity}
Let $X$ be a stratified $G$-space and let $E_G^*$ be a
multiplicative cohomology theory as above. Assume that the Euler
classes \(e(V_i) \in E_G^*(F_i)\) of the vector bundles \(V_i \to
F_i\) are not zero divisors. Then the inclusion \(\,\imath\!:
\coprod F_i \, \into \, X\) induces an injection
\begin{equation}\label{referee's 1}
\imath^*: E^*_G(X) \to \prod_i E^*_G(F_i).
\end{equation}
Moreover, let $E^*_G(X)$ be given the induced filtration under the
above inclusion. Then the associated graded $E_G^*$-module
\(QE^*_G(X)\) is isomorphic to (the direct product of) the ideals
generated by the Euler classes in the $E^*_G(F_i)$. Explicitly,
\begin{equation}\label{referee's 2}
QE^*_G(X) \cong \prod_i e(V_i) E^*_G(F_i).
\end{equation}
\end{theorem}

\begin{proof}
We first prove the theorem when the stratification of $X$ is
finite. This is done by induction on the length of the
stratification.

We first consider the assertion that \eqref{referee's 1} is
injective.  If the length of the stratification is 0, then $X$ is
empty, both sides of (\ref{referee's 1}) are zero, and the result
trivially holds. We now argue the inductive step. Assume that the
stratification of $X$ has length $i$ (i.e. $X=X_i$) and consider
the cofiber sequence
\begin{equation}\label{cofibration}
X_{i-1}\longrightarrow  X_i \stackrel{p}{\longrightarrow} Th(V_i).
\end{equation}
It follows from the  assumption on the Euler class that
the long exact sequence in $E$-cohomology associated to (\ref{cofibration}) splits into short exact sequences
\begin{equation}\label{SES}
0 \longrightarrow E_G^*(Th(V_i)) \stackrel{p^*}{\longrightarrow}
E_G^*(X_i) \longrightarrow E_G^*(X_{i-1}) \longrightarrow 0.
\end{equation}
To see this, we prove that  $p^*$ is an injection. Indeed, the
composition
\[
\xymatrix{
E^*_G(F_i)  \ar[r]^(0.4){\cdot u}_(0.4){\cong} &  E^*_G(Th(V_i)) \ar[r]^(0.6){p^*} &
E^*_G(X_i) \ar[r] & E^*_G(F_i)
}
\]
is multiplication by the Euler class $e(V_i)$, and is therefore
injective. The first map is the Thom isomorphism (see
\cite[Theorem~9.2]{May}), so the middle map $p^*$ must be
injective.

Now consider the map of short exact sequences
\begin{equation}\label{eq:shortExact}
\begin{array}{c}
\begin{xy}
\xymatrix{
0 \ar[r]& E_G^*(Th(V_i)) \ar[r]%\ar[d]
& E_G^*(X_i) \ar[r]%\ar[d]
& E_G^*(X_{i-1}) \ar[r]%\ar[d]
& 0 \\
0 \ar[r]& E_G^*(F_i) \ar[r]& \displaystyle\prod_{j \leq
i}^{\phantom{boo}} E_G^*(F_j) \ar[r]& \displaystyle\prod_{j <
i}^{\phantom{boo}} E_G^*(F_j) \ar[r]& 0. }
 \POS "1,2",\ar@{}"2,2",\ar"2,2"!E+/0pt/*{}
 \POS "1,3",\ar@{}"2,3",\ar"2,3"!E+/11pt/*{}
 \POS "1,4",\ar@{}"2,4",\ar"2,4"!E+/11pt/*{}
\end{xy}
\end{array}
\end{equation}
The left vertical map is injective by the assumption on $e(V_i)$,
with image $e(V_i)E_G^*(F_i)$. The right vertical map is injective
by induction. By the Five Lemma, the central map is also
injective. This proves \eqref{referee's 1} when the filtration of
$X$ is finite.

We now prove \eqref{referee's 2}.  Again, the base case is trivial,
since both sides of \eqref{referee's 2} are zero when the
stratification has length zero.  We now argue the inductive step. The
associated graded $QE^*_G(X_i)$ is isomorphic to $E^*_G(Th(V_i))\oplus
QE^*_G(X_{i-1})$. The image of $QE^*_G(X_{i-1})$ under the rightmost
vertical map in \eqref{eq:shortExact} is $\prod_{j<i}e(V_j)E^*_G(F_j)$
by the induction hypothesis. So, the image of $QE^*_G(X_i)$ under the
center vertical map is
$$
QE^*_G(X_i) \cong e(V_i)E^*_G(F_i)\oplus \prod_{j<i} e(V_j)E^*_G(F_j) = 
\prod_{j\le i} e(V_j)E^*_G(F_j),$$ 
as claimed in \eqref{referee's 2}.

For both statements \eqref{referee's 1} and \eqref{referee's 2}, the general case $\displaystyle
X=\underrightarrow{\lim} \, X_i$ follows directly from the finite
case since
\[
E^*_G(X)= \underleftarrow{\lim} \, E^*_G(X_i).
\]
Note that there is no Milnor \(\lim^1\) term here because the maps
$E_G^*(X_i)\to E_G^*(X_{i-1})$ are all surjective.
\end{proof}

\section{The combinatorial description of $E_G^*(X)$}\label{se:GKM}

We now identify the image of $E^*_G(X)$ in $\prod E^*_G(F_i)$: it
is specified by simple combinatorial restrictions. This is the
content of Theorem~\ref{th:GKM}.  In order to
make this computation, we must make some 
additional assumptions on $X$. 
We formalize our hypotheses on $X$ below.

\begin{assumption}\label{as:injectivity}
The space $X$ is equipped with a $G$-invariant
stratification
\[
X = \bigcup_{i\in I} X_i
\]
and each successive quotient \(X_i/X_{<i}\) is homeomorphic
  to the Thom space of a $G$-equivariant
  vector bundle \(\pi_i:V_i\to F_i\). Here $X_{<i}$ denotes the subspace $\bigcup_{j<i}X_j\subset X_i$.
\end{assumption}

\begin{assumption}\label{as:directSum}
The bundles $V_i \to F_i$ are $E$-orientable and admit $G$-equivariant direct sum
decompositions \[\big(\pi_i:V_i\to F_i\big) \cong \bigoplus_{j < i} \big(\pi_{ij}:V_{ij}\to F_i\big)\]
into $E$-orientable vector bundles $V_{ij}$. We allow the case $V_{ij}=0$.
\end{assumption}

\begin{assumption}\label{as:attachMap}
There exist $G$-equivariant maps \(f_{ij}: F_i \to F_j\) such that
the attaching maps $\varphi_i:S(V_i)\to X_{i-1}$, when restricted
to $S(V_{ij})$, are given by \[\varphi_i|_{S(V_{ij})}=f_{ij} \circ
\pi_{ij}.\] Here, we identify the $F_j$ with their images in
$X_{i-1}$.
\end{assumption}

\begin{assumption}\label{as:eulerClass}
The Euler classes \(e(V_{ij})\) are not zero divisors and are
pairwise relatively prime in $E^*_G(F_i)$. Namely, for any class
\(x \in E^*_G(F_i),\) we have that
\[(\forall j)\;e(V_{ij}) | x  \;\Leftrightarrow\; e(V_i) | x.\]
\end{assumption}

With these assumptions, we may now formulate our main theorem.

\begin{theorem}\label{th:GKM}
Let $X$ be a $G$-space satisfying Assumptions~\ref{as:injectivity} through
\ref{as:eulerClass}. Then the map
\[
\imath^*:  E_G^*(X) \to \prod_i E_G^*(F_i)
\]
is injective with image
\begin{equation}\label{eq:iotaImage}
R := \left.\left\{ (x_i) \in \prod_i E^*_G(F_i)\  \right|\  e(V_{ij}) \mid x_i -
f_{ij}^*(x_j) \ \text{for all}\ j<i \right\}.
\end{equation}
\end{theorem}

When $V_{ij}=0$ in the theorem above, the relation $e(V_{ij}) \mid
x_i - f_{ij}^*(x_j)$ is vacuous because $e(0)=1$. We introduce a
decorated graph $\Gamma$ that carries all the information from $X$
necessary to compute the image $R$ of $E_G^*(X)$.  Each edge
of $\Gamma$ corresponds to a non-vacuous relation.

\begin{definition}\label{GKMgraph}
The GKM graph $\Gamma$ associated to $X$ is the graph with
one vertex $v_i$ for each subspace $F_i$ and an edge $(v_i, v_j)$
whenever  $V_{ij}$ is non-zero. Each edge is labeled with the
bundle $V_{ij}$ and the map $f_{ij}:F_i\to F_j$.
\end{definition}

\begin{remark}
In Sections~\ref{sec:eqvtCohAffGr} and \ref{sec:examples}, the
description of $\Gamma$ simplifies greatly. In those examples, all
the $F_i$ are single points, and the maps $f_{ij}:F_i\to F_j$ are
the only possible ones. Moreover, the bundles $V_{ij}$ are all
1-dimensional complex $T$-representations. Hence $\Gamma$ is a
graph with a character $\alpha\in\Lambda:=\Hom(T,S^1)$ attached
to each edge.
\end{remark}

\begin{remark}
Theorem~\ref{th:GKM} generalizes many results found in the
literature.  We survey some of these results here.
\begin{enumerate}
\item[{A.}] Suppose that $X$ is a projective variety equipped with
an algebraic action of a complex torus, with finitely many $0$-
and $1$-dimensional orbits. Let $E_G^*$ be ordinary
$T$-equivariant cohomology.  In this setting, Theorem~\ref{th:GKM}
is precisely the result of Goresky, Kottwitz, and MacPherson
\cite{GKM}.

\item[{B.}] Theorem~\ref{th:GKM} recovers the main theorem of
\cite{GH04} when $X$ is a compact Hamiltonian $T$-space with
possibly non-isolated fixed points, and generalizes this result to
equivariant $K$-theory.

\item[{C.}] When $E_G^*$ is $T$-equivariant $K$-theory with
complex coefficients and $X$ is a GKM manifold, then
Theorem~\ref{th:GKM} is identical to
\cite[Corollary~A.5]{KnuRos03}.

\item[{D.}] If $X$ is a Kac-Moody flag variety and $E_G^*$ is
$T$-equivariant $K$-theory, then Theorem~\ref{th:GKM} is closely
related to a result of Kostant-Kumar \cite{KosKum87}. Indeed,
their Theorem~3.13 identifies $K_T^*(\G/\B)$ with the subring of
elements of $\prod_W K_T^*$ that are mapped to $K_T^*$ by certain
operators, which include the divided difference operators
$$(\delta_w -\delta_{w{r_\alpha}})\frac{1}{1-e^{\alpha}}$$ for all
$w\in W$ and reflections $r_\alpha$.  These are exactly the same
conditions as in \eqref{eq:iotaImage}.  Their Corollary~3.20
determines $K_T^*(\G/\P)$ in a similar fashion.
\end{enumerate}
\end{remark}

Before proving Theorem \ref{th:GKM}, we give a Lemma which
computes $E_G^*(X)$ when the stratification of $X$ has length 2.

\begin{lemma}\label{lem:FjVij}
Let $Y=F_1\cup_{\varphi}D(V)$ be obtained by gluing the sphere
bundle of $\pi:V\to F_2$ onto $F_1$, where $\varphi=f\circ\pi$ for
a map $f:F_2\to F_1$. Assume that $e(V)$ is not a zero divisor.
Then the images of the restriction maps \(\imath^*: E^*_G(Y, F_1)
\to E^*_G(F_2)\) and \(\jmath^*: E^*_G(Y) \to E^*_G(F_1) \oplus
E^*_G(F_2)\) are
\begin{equation}\label{i}
\imath^*(E^*_G(Y, F_1)) = \left\{ g \in E^*_G(F_2)\ \Big| \ e(V)
\mid g \right\}
\end{equation}
and
\begin{equation}\label{j}
\jmath^*(E^*_G(Y)) = \left\{ (g_1,g_2) \in E^*_G(F_1) \oplus
E^*_G(F_2)\ \Big| \ e(V) \mid g_2-f^*(g_1) \right\},
\end{equation}
respectively.
\end{lemma}

\begin{proof}
Clearly $E_G^*(Y,F_1) \cong E_G^*(Th(V)) \cong E_G^*(F_2)$ via the
Thom isomorphism. The map $$E_G^*(F_2)\cong
E_G^*(Th(V))\stackrel{\imath^*}{\longrightarrow}E_G^*(F_2)$$ is
multiplication by $e(V)$, so $\Im(\imath^*)$ is $e(V)E^*_G(F_2)$
as claimed in (\ref{i}).

The space $Y$ retracts onto $F_1$ via the map $f\circ\pi$, so the long exact sequence associated to the pair $(Y,F_1)$ splits.
Now consider the diagram
\[
\xymatrix{
0\ar[r]& E^*_G(Y,F_1)\ar[r]\ar[d]^{\imath^*}& E^*_G(Y)\ar[r]\ar[d]^{\jmath^*} & E^*_G(F_1)\ar[r]\ar[d]& 0 \\
0\ar[r]& E^*_G(F_2)  \ar[r]      & E^*_G(F_1\sqcup F_2)\ar[r]&
E^*_G(F_1)\ar[r]      & 0. }
\]
Both rows split, and we get $\Im(\jmath^*)=E^*_G(F_1)\oplus
\Im(\imath^*)$, where $E^*_G(F_1)$ is mapped via the diagonal
inclusion $(1,f^*): E_G^*(F_1) \to E^*_G(F_1)\oplus E^*_G(F_2)$.
It is now straightforward to check that
$\{(g_1,f^*(g_1))\}\oplus\{(0,g_2):e(V)\ |\ g_2\}$ is the same
group as described in (\ref{j}).
\end{proof}

We now have the technical tool to prove our main theorem.
\medskip

\noindent{\bf Proof of Theorem \ref{th:GKM}:} The map $\imath^*$
is injective by Theorem~\ref{thm:injectivity}, so we must show
that its image $\Im(\imath^*)$ equals the ring $R$ of
\eqref{eq:iotaImage}.

We first show that $\Im(\imath^*)\subseteq R$. Let $Y_{ij}$ be the subspace of $X$ given by
\[
Y_{ij}:= F_j\cup_{f_{ij}\circ\pi_{ij}}D(V_{ij}).
\]
Consider a class $x\in E^*_G(X)$, and let $x_i$ denote its
restriction to $F_i$. Since $(x_j,x_i)$ is the image of
$x|_{Y_{ij}}\in E^*_G(Y_{ij})$ under the restriction map
$E^*_G(Y_{ij})\to E^*_G(F_j)\oplus E^*_G(F_i)$, we know by
Lemma~\ref{lem:FjVij} that
\begin{equation}\label{kappaij}
e(V_{ij}) \mid x_i - f_{ij}^*(x_j).
\end{equation}
The conditions (\ref{kappaij}) characterize $R$, so we conclude
$(x_i)\in R$.

We now have a map $E^*_G(X)\to R$ and want to show that it is
surjective. Following Remark~\ref{re:posets}, we are using $I=\N$.
We argue by induction on the length of the stratification. If the
length is zero, then $X=\emptyset$ and there is nothing to show.
We now assume that $X=X_i$ and that surjectivity holds for
$$
E^*_G(X_j)\to R_j:=\left.\left\{ (x_k) \in \prod_{k\le j}
E^*_G(F_k)\ \right|\  e(V_{k\ell}) \mid x_k - f_{k\ell}^*(x_\ell)
\ \text{for all}\ \ell<k \right\}
$$
for all $j<i$.

Let \(r_i: R_i \to R_{i-1}\) be the restriction map. By Assumption 
\ref{as:eulerClass}, its kernel can be written
\begin{equation}\label{eq:kerr}
\mathrm{ker}(r_i) = \left.\left\{ (x_j) \in \prod_{j \leq i}
E^*_G(F_j)\ \right| \ \begin{tabular}{l} \(x_j = 0\) for \(j < i\) \\
\(e(V_{ij}) \mid x_i\) for all \(j < i\) \\
\end{tabular} \right\}\simeq e(V_i)E^*_G(F_i).
\end{equation}
We now consider the following commutative diagram:
\begin{equation}\label{eq:FiveLemmadiag}
\begin{array}{c}
\xymatrix{
0 \ar[r] & E^*_G(X_{i}, X_{i-1}) \ar[r] \ar[d] & E^*_G(X_{i})
\ar[r]\ar[d] & E^*_G(X_{i-1}) \ar[r] \ar[d] & 0 \\
0 \ar[r] & \mathrm{ker}(r_i) \ar[r] &
R_i \ar[r]^{r_i} & R_{i-1}. \\
}
\end{array}
\end{equation}
The top sequence comes from the long exact sequence of the
pair, which splits into short exact sequences as shown in the
proof of Theorem \ref{thm:injectivity}. By the induction
hypothesis, we know that the right vertical arrow is an
isomorphism. By comparing (\ref{i}) and (\ref{eq:kerr}), the left
vertical arrow is also an isomorphism. It is now an easy diagram
chase to verify that $r_i$ is surjective and that
$E^*_G(X_i)\simeq R_i$.

Finally, we note that
$$
E^*_G(X) = \underleftarrow{\lim}\ E^*_G(X_i) =
\underleftarrow{\lim}\ R_i = R,
$$
completing the proof.
\qed

\section{Module generators}\label{se:generators}

The second part of Theorem \ref{thm:injectivity} gives us a lot of
information about the structure of $E^*_G(X)$ as an
$E^*_G$-module. When the spaces $F_i$ consist of isolated fixed
points, we can say more. With this assumption, (\ref{referee's 2})
tells us that as an $E^*_G$-module, $E^*_G(X)$ is
(non-canonically) a product of principal ideals of $E^*_G$:
\[
E^*_G(X)\cong \prod_{v\in F} e(V_v)E^*_G,
\]
where $F=\cup F_i$ and $V_v$ is the fiber over $v$. Moreover,
given a collection of classes $x_v\in E^*_G(X)$, one for each
$v\in F$, it is very easy to check whether they form a set of free
generators\footnote{Here $E^*_G(X)$ should be viewed as a
topological $E^*_G$-module, and the word `generator' should be
interpreted in the topological sense.} for $E^*_G(X)$. 

We write $v<w$ when $v\in F_i$, $w\in F_j$ and $i<j$. 
We write $v\le w$ if
$v<w$ or $v=w$.  Let $x_{v}(w)$ denote $x_v|_{w}$. We then have:
\begin{proposition}\label{gen1}
Suppose $X$ satisfies Assumptions 1-4 and that the spaces $F_i$
consist of isolated fixed points. Let $x_v\in E^*_G(X)$ be classes
satisfying
\begin{equation}\label{xij}
\begin{array}{l}
x_{v}(w)=0 \ \text{for}\ w\not\ge v;\ \mbox{and}\\
x_{v}(v) \;\; \text{is a generator of the ideal} \;\;e(V_v)E^*_G.
\end{array}
\end{equation}
Then $\{x_v\}$ is a set of free topological $E^*_G$-module
generators.\qed
\end{proposition}

It might happen that a space $X$ with $G$-action satisfies the Assumptions 1-4 for
some cohomology theory $E^*_G$, but that Assumption
\ref{as:eulerClass} fails for some closely related cohomology
theory $\widetilde{E}^*_G$. For example, this can happen when
$\widetilde{E}^*_G$ is non-equivariant $E$-cohomology $E^*(X):=
E^*_G(X\times G)$, or when $E^*_G=H^*_G(-;\Z)$ and
$\widetilde{E}^*_G=H^*_G(-;\Z/2)$. In that case we have:

\begin{proposition}\label{modulespectrum}
Suppose $X$ satisfies Assumptions 1-4 for the cohomology theory
$E^*_G$, and that the $F_i$ consist of isolated fixed points. Let
$\widetilde{E}^*_G$ be a module cohomology theory over the ring
cohomology theory $E^*_G$. Then one can recover
$\widetilde{E}^*_G(X)$ by tensoring
\[
\widetilde{E}^*_G(X)=E^*_G(X)\widehat\otimes_{E^*_G}\widetilde{E}^*_G.
\]
Here $E^*_G(X)$ is viewed as a topological $E^*_G$-module and
$\widehat\otimes$ denotes the completed tensor product.

In particular, if $\widetilde{E}^*_G$ is an $E^*_G$-algebra and
$x_v\in E^*_G(X)$ satisfy (\ref{xij}), then $x_v\otimes 1$ are
free $\widetilde{E}^*_G$-module generators of
$\widetilde{E}^*_G(X)$.
\end{proposition}

\begin{proof}
We argue by induction on the length of the stratification.

Without loss of generality, we may assume the $F_i$ are single
points. The short exact sequence (\ref{SES}) consists of free
$E^*_G$-modules. Therefore, the functor
$-\otimes_{E^*_G}\widetilde{E}^*_G$ preserves exactness, and we
get the following commutative diagram:
\[
\xymatrix{ 0\ar[r]& E^*_G(Th(V_i))\otimes_{E^*_G}
\widetilde{E}^*_G\ar[r]\ar[d]& E^*_G(X_i)\otimes_{E^*_G}
\widetilde{E}^*_G\ar[r]\ar[d]&
E^*_G(X_{i-1})\otimes_{E^*_G} \widetilde{E}^*_G\ar[r]\ar[d] & 0\\
&\widetilde{E}^*_G(Th(V_i))\ar[r]^{\alpha}&\widetilde{E}^*_G(X_i)\ar[r]^{\beta}&\widetilde{E}^*_G(X_{i-1}).&
}
\]
The right vertical arrow is an isomorphism by induction. The left
vertical arrow is an isomorphism since
\[
E^*_G(Th(V_i))\otimes_{E^*_G} \widetilde{E}^*_G\cong
E^*_G(F_i)\otimes_{E^*_G} \widetilde{E}^*_G\cong
\widetilde{E}^*_G\cong \widetilde{E}^*_G(Th(V_i)),
\]
where the first and last isomorphisms are the equivariant suspension isomorphisms.

A diagram chase shows that $\beta$ is surjective, so the bottom
long exact sequence splits and the map $\alpha$ is injective. We
deduce by the Five Lemma that the middle vertical map is also an
isomorphism, as desired.

Finally, if the filtration is infinite, we have
\begin{align*}
\widetilde{E}^*_G(X)=\underleftarrow\lim \,\widetilde{E}^*_G(X_i)=&\underleftarrow\lim \big(E^*_G(X_i)\otimes_{E^*_G}\widetilde{E}^*_G\big) \\
=&\big(\underleftarrow\lim
\, E^*_G(X_i)\big)\widehat\otimes_{E^*_G}\widetilde{E}^*_G=E^*_G(X)\widehat\otimes_{E^*_G}\widetilde{E}^*_G.
\end{align*}
\vskip -0.3in
\end{proof}

Assume now that $X$ is a CW complex with $G$-invariant cells\footnote{Careful: we don't mean tht $X$ is a $G$-CW complex.}, that the filtration \eqref{stratification} is the
usual filtration by skeleta (indexed by $\N$), and that
$E^*_G(X)=H^*_G(X):=H^*(X\times_G EG)$ is ordinary equivariant
cohomology.  In this case, we can give a canonical set of free
generators for $H^*_G(X)$. As before, we let $F=\cup F_i$, where $F_i$ 
is now the set of the centers of the
$i$-dimensional cells. We write
$|v|=i$ whenever $v\in F_i$ and recall the notation $x_v(w)$ for $x_v|_w$.

\begin{proposition}\label{gen3}
Let $X$ be a CW complex as above. Then there is a unique set
$\{x_v\}_{v\in F}$ of free generators for the $H_G^*$-module
$H_G^*(X)$ satisfying the conditions:

\begin{enumerate}
\item each $x_v$ is homogeneous of degree $|v|$; \item if $|w|\le
|v|$, $w\not= v$, then \(x_{v}(w)= 0 \in H_G^*;\) and
\item the element \(x_{v}(v)\) is the equivariant Euler class
$e(V_v):=e(V_v\times_G EG\to BG)\in H^*_G$, where $V_v$ is the
cell of $X$ with center $v$.
\end{enumerate}
\end{proposition}

\begin{proof}
We first construct the classes $x_v$. Assume by induction that we
have classes $x'_w$ in $H^*_G(X_{i-1})$ for  $|w|<i$. To extend
these to $H^*_G(X_i)$, consider the short exact sequence
\[
\xymatrix{ 0 \ar[r] & H^*_G(X_i, X_{i-1}) \ar[r] & H^*_G(X_i)
\ar[r] & H^*_G(X_{i-1}) \ar[r] & 0 }
\]
and note that \[H^*_G\big(X_i, X_{i-1}\big) \cong
H^*_G\Big(\bigvee_{|v|=i} Th(V_v)\Big) \cong \prod_{|v|=i}
H^*_G\big(Th(V_v)\big).\] The spaces $Th(V_i)$ are $G$-spheres,
so each \(H^*_G(Th(V_v))\) has a canonical generator $u_v$. The
restriction of $u_v$ to the center $v$ of $V_v$ is the equivariant
Euler class $e(V_v)$. The classes $x'_w$ of $H^{*}_G(X_{i-1})$
have a unique lift $x_w$ to $H^{*}_G(X_i)$ because $H^{k}_G(X_{i},
X_{i-1})$ is zero for all $k<i$. It is straightforward to check
that these lifts, along with the images $x_v$ of the chosen
generators $u_v$ of $H^{*}_G(X_{i}, X_{i-1})$, satisfy the above
conditions and generate $H^*_T(X_i)$. We take a limit over $i$ to obtain the
generators $x_v\in H_T^*(X)$.

We show that Conditions 1, 2 and 3 characterize the generators
$x_v$. Let $\{\widetilde{x}_v\}$ be another set of generators
satisfying the same conditions. Write them as \(\widetilde{x}_v =
\sum_w b_{vw} x_w.\) By Condition 2, we have \(b_{vw} = 0\)
whenever $|w| \leq |v|$ and $w \neq v$. By Condition 3, \(b_{vv} =
1\). Finally, \(b_{vw} = 0\) when \(|w| > |v|,\) because otherwise
$\widetilde{x}_v$ would not be homogeneous.
\end{proof}
\begin{remark}
Suppose $X$ is a manifold with a $G$-invariant Morse function $f$
and a CW decomposition constructed from the Morse flow. Then the
above construction is the same as the following: given a fixed
point $v$, consider the flow-up manifold $\Sigma_v$ of codimension
\(|v|\). By Poincar\'e duality, it represents a cohomology class
$x_v$. It is straightforward to see that the $x_v$ satisfy
Conditions~1, 2 and 3 of Proposition \ref{gen3}.
\end{remark}

\begin{remark}\label{rem:KthGen}
There are other situations when it is possible to find canonical
module generators. For example, such generators exist when $X$ is
a complex algebraic variety or a symplectic manifold, and $E_G^*$
is equivariant $K$-theory.  The algebraic construction involves
resolving the structure sheaf of the ``flow-up'' varieties
$\Sigma_v$. See \cite{BFM} for details.  The symplectic
construction can be found in \cite{GK03}.
\end{remark}

We illustrate these generators for some examples in
Section~\ref{sec:examples}.

\section{Kac-Moody flag varieties}\label{sec:eqvtCohAffGr}

We now turn our attention to the main examples that motivate the
results in this paper. These are homogeneous spaces $\G/\P$ for a
(not necessarily symmetrizable) Kac-Moody group $\G$, defined over
$\C$, with $\P$ a parabolic subgroup.  Specific examples of such
homogeneous spaces include finite dimensional Grassmannians, flag
manifolds, and based loop spaces $\Omega K$ of compact simply
connected Lie groups $K$.

We first take a moment to explicitly describe $\Omega K$ as a
homogeneous space $\G/\P$. Let $LK$ be the group of polynomial
loops
\[
LK := \{ \gamma:S^1\to K \},
\]
where the group structure is given by pointwise multiplication. By
polynomial, we mean that the loop is the restriction $S^1 = \{z
\in \C: |z|=1\} \to K$ of an algebraic map \(\C^* \to K_{\C}.\)
The space of based polynomial loops is defined by
\[
\Omega K := \{ \sigma \in LK |\ \sigma(1) = 1\in K \}.
\]
The group $LK$ acts transitively on $\Omega K$ by
\begin{equation}\label{eq:LSutonPSut}
(\gamma \cdot \sigma)(z) = \gamma(z) \sigma(z) \gamma(1)^{-1}.
\end{equation}
The stabilizer of the constant identity loop is exactly $K$, the
subgroup of constant loops. Thus $\Omega K\cong LK/K$.

Now let $\G$ be the affine Kac-Moody group $\G = \widehat{LK_{\C}}
\rtimes \C^*$. Here, $LK_{\C}$ is the group of algebraic maps
\(\C^* \to K_{\C}\), \(\widehat{LK_{\C}}\) is the universal
central extension of $LK_{\C}$, and the $\C^*$ acts on $LK_{\C}$
by rotating the loop. The parabolic $\P$ is $\widehat{L^+K_{\C}}
\rtimes \C^*$, where $L^+K_{\C}$ is the subgroup of $LK_{\C}$
consisting of maps \(\C^* \to K_{\C}\) that extend to maps \(\C
\to K_{\C}.\) It is shown in \cite[8.3]{PS} that $\Omega K$ can be
identified as a homogeneous space $\G/\P$. We briefly sketch this
argument. The group $LK$ acts on $\G/\P$ by left multiplication,
and the stabilizer of the identity is $\P \cap LK$. This
intersection is the set of polynomial maps \(\C^* \to K_{\C}\)
which extend over $0$, and which send $S^1$ to $K$. Thus, a loop
\(\gamma\) in $\P \cap LK$ satisfies  the condition \(\gamma(z) =
\theta(\gamma(1/\bar{z})),\) where $\theta$ is the Cartan
involution on $K_{\C}$. Therefore, since $\gamma$ extends over
zero, by setting \(\gamma(\infty) = \theta(\gamma(0)),\) it also
extends over $\infty$. But then $\gamma$ is an algebraic map from
\(\PP^1\) to \(K_{\C},\) and is therefore constant, since $K_{\C}$
is affine. Hence \(\P \cap LK = K.\)

\begin{remark}
We have only considered the space of polynomial loops in $K$.
However, our results still apply to other spaces of loops, such as
smooth loops, 1/2-Sobolev loops, etc.  Indeed, the polynomial
loops are dense in these other spaces of loops \cite[3.5.3]{PS},
\cite{Mitch}. By Palais' theorem \cite[Theorem 12]{Pal}, these
dense inclusions are weak homotopy equivalences. The inclusions of
$T'$-fixed point sets for $T'$ a closed subgroup of $T$ are also
equivalences.  So the various forms of $\Omega K$ are actually
equivariantly weakly homotopy equivalent.
\end{remark}

Let us return to the general case. Let $T_\G$ be the maximal torus
of $\G$. The center $Z(\G)$ acts trivially on $X=\G/\P$, so the
quotient group $T:=T_\G/Z(\G)$ acts on $X$. We need to check that
this space $X$ with this $T$-action satisfy
Assumptions~\ref{as:injectivity}-\ref{as:eulerClass} that are the
hypotheses of Theorem~\ref{th:GKM}.  It is known (see for
example \cite{BilDye, KP, KosKum87, Mitch}) that $\G/\P$ admits a
$T$-invariant CW decomposition
\begin{equation}\label{eq:GPdecomp}
\G/\P = \coprod_{[w] \in W_\G/W_\P} \B\tilde{w}\P/\P,
\end{equation}
where $W_\G$ and $W_\P$ are the Weyl groups of $\G$ and of (the
semisimple part of) $\P$ respectively, and $\tilde{w}$ is a
representative of $w$ in $\G$. This is the filtration of
Assumption~\ref{as:injectivity}. Each cell is homeomorphic to a
$T$-representation and has a single $T$-fixed point \(\bar{w} :=
\tilde{w}\P/\P\) at its center. These cells are the $V_i$ and the
fixed points are the $F_i$. The $T$-representation $V_i$ is
isomorphic to the tangent space
\[
\mathrm{T}_{\bar{w}}\B\bar{w} =
\mathrm{T}_{\bar{w}}\B\tilde{w}\P/\P = \mathfrak{b}/\mathfrak{b}
\cap \tilde{w}\mathfrak{p}\tilde{w}^{-1} =
\mathfrak{b}/\mathfrak{b} \cap w \cdot \mathfrak{p}.
\]
This tangent space decomposes into 1-dimensional representations,
corresponding to the roots contained in $\mathfrak{b}$ but not in
$w \cdot \mathfrak{p}$.  These subspaces are the $V_{ij}$ of
Assumption~\ref{as:directSum}.

We now check Assumption~\ref{as:attachMap}. Since the $F_i$ are
points, we only need to show that the attaching map $\varphi_i:
S(V_i) \to X_{i-1}$ maps each $S(V_{ij})$ onto the point $F_j$. In
other words, we need to show that the closure of $V_{ij}$ is a
$2$-sphere with north and south poles $F_i$ and $F_j$. Pick a root
$\alpha$ in $\mathfrak{b}$ but not in $w \cdot \mathfrak{p}$. Let
$e_{\alpha}, e_{-\alpha}$ be the standard root vectors for
$\alpha, -\alpha$. Let \(SL(2,\C)_{\alpha}\) be the subgroup of
$\G$ with Lie algebra spanned by $e_{\alpha}, e_{-\alpha}$ and
$[e_{\alpha},e_{-\alpha}]$, and let $\B_{\alpha}$ be the Borel of
$SL(2,\C)_{\alpha}$ with Lie algebra spanned by $e_{\alpha}$ and
$[e_{\alpha},e_{-\alpha}]$. Let \(\tilde{r}_{\alpha} :=
exp(\pi(e_{\alpha} - e_{-\alpha})/2)\) represent the element
$r_{\alpha}$ of the Weyl group which is reflection along $\alpha$.
Let $F_i$ be the point $\bar{w}$ and $F_j$ the point
$r_\alpha\bar{w}$. The $\alpha$-eigenspace in the cell $\B\bar{w}$
is $\B_{\alpha}\bar{w} = V_{ij} \cong \C$. Its closure is
$SL(2,\C)_{\alpha} \bar{w} \cong \PP^1$, and the point at infinity
is given by $\tilde{r}_{\alpha} w\P/\P = r_{\alpha} \bar{w} =
F_j$, as desired.

Finally, we need to check Assumption~\ref{as:eulerClass}.  To do
this, we must show that for the roots contained in $\mathfrak{b}$
but not in $w \cdot \mathfrak{p}$, the corresponding Euler classes
are pairwise relatively prime. This is true for a large class of
$T$-equivariant complex oriented cohomology theories including
$H_T^*(-;\Z)$, $K_T^*$ and $MU_T^*$.

\begin{lemma}\label{le:Ass4forG/P}
Let $E_T^*$ be $H_T^*(-;\Z)$, $K_T^*$ or $MU_T^*$. Let $\alpha_i$
be any finite set of non-zero characters such that no two are
collinear. Moreover, if $E_T^*=H_T^*(-;\Z)$, assume that no prime
$p$ divides two of the $\alpha_i$. Then the corresponding Euler
classes $e(\alpha_i)$ are pairwise relatively prime in $E_T^*$.
\end{lemma}

\begin{proof}
The equivariant cohomology ring $H_T^*$ is the symmetric
algebra\footnote{This is true if one restricts the $RO(T)$-grading
of \cite{May} to the more familiar $\Z$-grading. Otherwise, one
has various periodicities with respect to all zero-dimensional
virtual $T$-representations.} $Sym^*(\Lambda)$ on the weight lattice
of $T$.  This is a unique factorization domain, and the Euler
classes $e(\alpha_i) = \alpha_i$ decompose into an integer times a
primitive character. The result follows immediately in this case.

The equivariant $K$-theory ring $K_T^0$ is the group ring
$\Z[\Lambda]$ generated by symbols $e^\alpha$.  For each $\alpha$
in our set of characters, let $\bar{\alpha}$ be the primitive
character in that direction, so $\alpha = n \bar{\alpha}$. The
Euler classes $e(\alpha_i) = 1-e^{\alpha_i}$ factorize as a
product of cyclotomic polynomials
$$
1 - e^{\alpha_i} =\ \prod_{d | n_i} \Phi_d(e^{\bar{\alpha}_i}).
$$
The factors $\Phi_d(e^{\bar{\alpha}_i})$ are all distinct, so the
result follows.

To prove the result about complex cobordism, we argue by induction
on the number of characters in our set. The base case is trivial.
Assume by induction that the result holds for $n$ characters and
that we are given a set $\alpha,\beta_1,\dots,\beta_n$ of $n+1$
characters satisfying the hypotheses of the lemma. Let $x$ be a
class in $MU_T^*$ which is divisible by each of the Euler classes
of the above characters. By induction, $x$ is divisible by the
product $\prod_i e(\beta_i)$, so there exists a class $b$ such
that $b \cdot \prod_i e(\beta_i)=x$.
We now consider the short exact sequence \cite[Theorem~1.2]{Sin01}
$$
\begin{xy}
\xymatrix{ 0 \ar[r] & MU_T^* \ar[r]^{\cdot e(\alpha)} & MU_T^*
%\ar[r]@{}^(0.4){res}
 & {}^{\phantom{|}}{MU_{\mathrm{Ker}(\alpha)}^*} \ar[r] & 0.}
\POS "1,3",\ar@{}"1,4",\ar"1,4"!E+/3pt/*{}^(0.65){res}
\end{xy}
$$

\noindent Since $x$ is divisible by
$e(\alpha)$, $$res(b)\cdot\prod_i res(\beta_i)=res(x) = 0.$$ By assumption, the restrictions
$\beta_i|_{\Ker(\alpha)}$ are non-torsion in the group of characters of $\Ker(\alpha)$. So by a result of
Sinha \cite[Theorem~5.1]{Sin01}
their Euler classes $e(\beta_i|_{\Ker(\alpha)})=res(e(\beta_i))$ are not zero divisors. We conclude that $res(b) = 0$.
Hence $b$ is a multiple of $e(\alpha)$, completing the proof.
\end{proof}

\begin{remark}
It is shown in \cite{CGK02} that any complex oriented
$T$-equivariant cohomology theory $E_T^*$ is an algebra over
$MU_T^*$. Combining this with Proposition~\ref{modulespectrum} and
Lemma~\ref{le:Ass4forG/P}, we may use our main
Theorem~\ref{th:GKM} to compute $E_T^*(\G/\P) =
MU_T^*(\G/\P) \widehat{\otimes}_{MU_T^*} E_T^*$.
\end{remark}

We conclude this section with an explanation of how to obtain the
pictures that we draw in Section~\ref{sec:examples}. The GKM graph
associated to $\G/\P$ has vertices $W_\G/W_\P$, with an edge
connecting $[w]$ and $[r_{\alpha}w]$ for all reflections
$r_{\alpha}$ in $W_\G$. The weight label on such an edge is
$\alpha$. It turns out that it is possible to embed this GKM graph
in $\mathfrak{t}^*$, the dual of the Lie algebra of $T$. Under
this embedding, the weight $\alpha_{ij}$ is then the primitive
element of $\Lambda\subset\mathfrak{t}^*$ in the direction of the
corresponding edge. To produce this embedding, we pick a point in
$\mathfrak{t}_\G^*$ whose $W_\G$-stabilizer is exactly $W_\P$,
take its $W_\G$-orbit, and draw an edge connecting any two
vertices related by a reflection in $W_\G$. This graph sits in a
fixed level of $\mathfrak{t}_\G^*$ (this is only relevant when
$\G$ is of affine type) and can therefore be thought of as sitting
in $\mathfrak{t}^*$.

\begin{figure}[h]
\begin{center}
\epsfig{figure=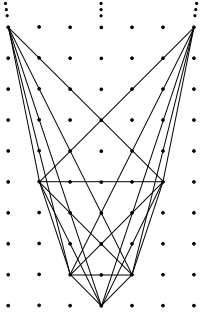,height=1.6in}
\end{center}
\begin{center}
\parbox{5in}{\caption{\small This is the GKM graph embedded in $\algt^*$ for
$\Omega SU(2)$, a homogeneous space for the loop group
$LSL(2,\C)$. }\label{fig:PSut}}
\end{center}
\end{figure}

These ideas are borrowed from the theory of moment maps in
symplectic geometry. In that context, $X$ is a symplectic manifold
with $T$-action and admits a moment map $\mu:X\to\mathfrak{t}^*$.
Consider the set $X^{(1)}$ of points with stabilizer of
codimension at most 1. The GKM graph is the image of $X^{(1)}$
under the moment map $\mu$. In our situation, $X^{(1)}$
corresponds exactly to the union of the $V_{ij}$.
Figure~\ref{fig:PSut} shows the image of the moment map for the
example $\Omega SU(2)$.

\section{Examples}\label{sec:examples}

\subsection{A homogeneous space for $G_2$}\label{sec:G2}

The complex Lie group $G_2$ contains two conjugacy classes of
maximal parabolic subgroups.  They correspond to the two simple
roots of $G_2$. We consider the case $X=G_2/\P$ and its natural
torus action, where $\P=\P_{long}$ is the parabolic generated by
the Borel subgroup and the exponential of the negative long simple
root. Equivalently, $X$ is the quotient of the compact group $G_2$
by a subgroup isomorphic to $U(2)$. The GKM graph is a complete
graph on 6 vertices and is embedded in $\mathfrak t^*\cong\R^2$
as a regular hexagon.

\begin{figure}[h]
%\psfrag{1}{$1$}
\psfrag{a}{$x$} \psfrag{b}{$y$} \psfrag{c}{$z$} \psfrag{d}{$s$}
\psfrag{e}{$t$}
\begin{center}
\epsfig{figure=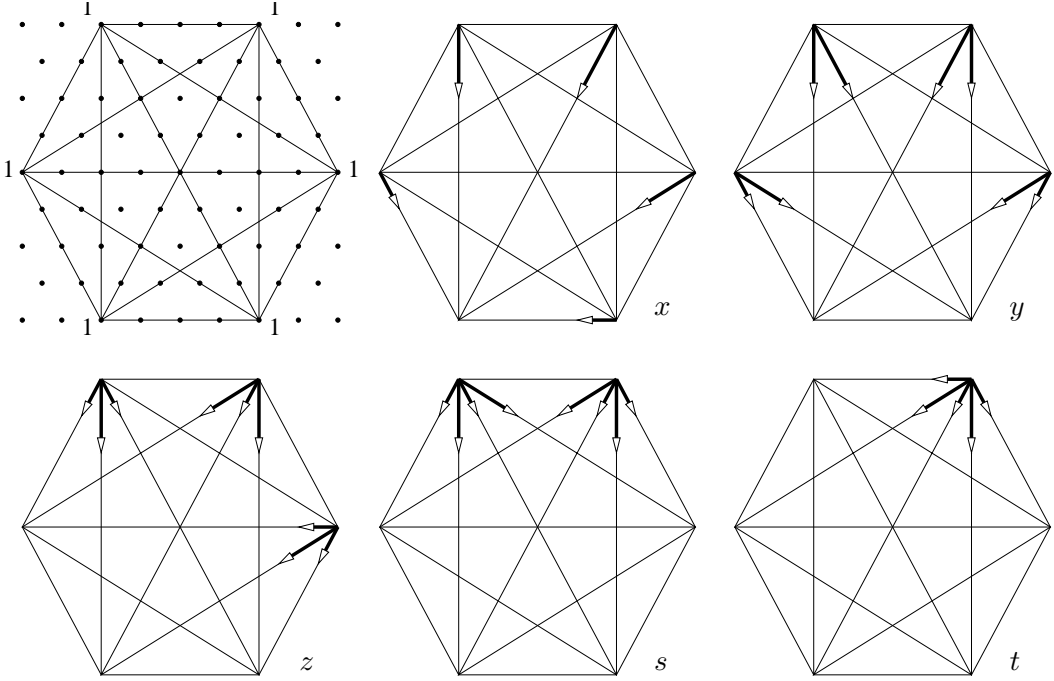,width=5.5in}
\end{center}
\begin{center}
\parbox{5in}{
\caption{ The module generators for $E_T^*(G_2/\P)$. We include
the lattice $\Lambda$ in the first diagram.}\label{fig:G2} }
\end{center}
\end{figure}

We now compute explicitly module generators $x_v$ of $E^*_T(X)$
for a large class of cohomology theories $E_T^*$, following
Section~\ref{se:generators}. We will represent them by their
restrictions $x_{v}(w):=x_v|_w$ to the various $T$-fixed points
$w\in F$. In this example, all the $x_{v}(w)$ happen to be Euler
classes of complex $T$-representations. This allows us to use the
following convenient notation to represent the classes $x_v$.  On
every vertex $w$ of $\Gamma$ we draw a bouquet of arrows
$\beta_j\in\Lambda$. By this, we mean that the class $x_{v}(w)\in
E^*_T(\{ w\})$ is the Euler class
\[
x_{v}(w)=e\big(\bigoplus_j\beta_j\big)=\prod_je(\beta_j).
\]
The vertices with no arrows coming out of them carry the class
$0$. Using these conventions, we draw the six module generators
$1,x,y,z,s,t$ of $E^*_T(G_2/\P)$ in Figure~\ref{fig:G2}.

Recall that Assumptions 1-4 are satisfied for the cohomology
theories $H^*_T(-;\Z)$, $K^*_T$ and $MU^*_T$, as shown in
Section~\ref{sec:eqvtCohAffGr}. To check that the elements shown
in Figure \ref{fig:G2} are module generators, we need to check two
things. First, we notice that the conditions (\ref{xij}) are
satisfied. Second, we need to verify that the elements $x,y,z,s,t$
satisfy the criteria \eqref{eq:iotaImage} for being elements of
$E^*_T(X)$.

To check \eqref{eq:iotaImage}, note that $e(\alpha)\in E^*_T$
divides $e(\beta)-e(\gamma)$ whenever $\beta-\gamma$ is a multiple
of $\alpha$  in $\Lambda$.  This is a trivial fact when $E^*_T$ is
ordinary $T$-equivariant cohomology or $T$-equivariant $K$-theory,
and is a consequence of the theory of equivariant formal group
laws when $E_T^*$ is an arbitrary $T$-equivariant complex oriented
cohomology theory \cite[p.\ 374]{CGK00}. Similarly $e(\alpha)$
divides a difference of products $\prod e(\beta_j)-\prod
e(\gamma_j)$ if the $\beta_j-\gamma_j$ are all multiples of
$\alpha$. Now, for each of the classes in Figure~\ref{fig:G2}, and
for each edge $(v,w)$ of $\Gamma$ with direction $\alpha$, we note
that the two bouquets of arrows $\{\beta_j\}$ at $v$ and
$\{\gamma_j\}$ at $w$ can be ordered in such a way that the
differences $\beta_j-\gamma_j$ are each in the direction of
$\alpha$. So we have checked (\ref{eq:iotaImage}) and hence by
Theorem \ref{th:GKM}, the classes in Figure~\ref{fig:G2} are
elements of $E^*_T(X)$. Thus, by Proposition~\ref{gen1}, they are
free module generators.

Even though the module generators look very similar in all
cohomology theories, the ring structures are different.  We
compute the ordinary $T$-equivariant cohomology and $K$-theory of
$X=G_2/\P$ to exhibit this phenomenon.

For cohomology theories $\widetilde{E}^*_T$ such as
$H^*_T(-;\Z/2)$, $H^*(-;\Z)$, $K^*$, or $MU^*$ for which
Assumption \ref{as:eulerClass} fails, we still have a good
understanding of $\widetilde{E}^*_T(X)$ by
Proposition~\ref{modulespectrum}.  We exploit this to compute
$H^*(X;\Z)$ from $H_T^*(X;\Z)$ and $K^*(X)$ from $K_T^*(X)$ below.

For the computation of $H^*_T(X;\Z)$, it is convenient to let
$a:=e(\put(0,3){\vector(1,0){12}}\hspace{4mm}),b:=e(\put(9,-2){\vector(-2,3){7}}\hspace{4mm})\in
H^2_T$ be the Euler classes of the characters
$\put(0,3){\vector(1,0){12}}\hspace{4mm},
\put(9,-2){\vector(-2,3){7}}\hspace{4mm}\in\Lambda$.  One then has
$H^*_T=\Z[a,b]$.  Using the embedding (\ref{referee's 1})
$H^*_T(X;\Z)\hookrightarrow\prod_F H^*_T$, we compute:
\begin{align*}
&x(x+a)=y, \\
&x(x+a)(x+b)=2z, \\
&x(x+a)(x+b)(x+2a+b)=2s,\mbox{ and}\\
&x(x+a)(x+b)(x+2a+b)(x+2b+a)=2t.
\end{align*}
To get the non-equivariant cohomology $H^*(X;\Z)$, it suffices by
Proposition~\ref{modulespectrum} to set $a=b=0$:
\begin{equation}\label{coh}
x^2=y, \quad
x^3=2z, \quad
x^4=2s, \quad
x^5=2t, \quad
x^6=0.
\end{equation}

In $K$-theory, it is more convenient to let $a,b\in K^0_T$ be the
characters $\put(0,3){\vector(1,0){12}}\hspace{4mm}$ and
$\put(9,-2){\vector(-2,3){7}}\hspace{4mm}\in\Lambda$ themselves
(not their Euler classes). We then have $K^0_T=\Z[a,a^{-1},b,b^{-1}]$, and all
other $K$-groups are either zero or isomorphic to $K^0$. We use
the convention that the Euler class of a line bundle $L$ is $1-L$.
We can now compute:
\begin{align*}
&x(ax+1-a)=y, \\
&x(ax+1-a)(bx+1-b)=(1+a^{-1})z-a^{-1}s, \\
&x(ax+1-a)(bx+1-b)(a^2bx+1-a^2b)=(1+b^{-1})s-b^{-1}t,\mbox{ and} \\
&x(ax+1-a)(bx+1-b)(a^2bx+1-a^2b)(ab^2x+1-ab^2)=(1+a^{-1}b^{-1})t.
\end{align*}
To get the non-equivariant $K$-theory, we set $a=b=1$ according to
Proposition~\ref{modulespectrum}: 
\begin{equation}\label{kot}
x^2=y, \quad
x^3=2z-s, \quad
x^4=2s-t, \quad
x^5=2t, \quad
x^6=0.
\end{equation}

\noindent We note that, as expected, the cohomology ring
(\ref{coh}) of $G_2/\P$ is the associated graded of the $K$-theory
ring (\ref{kot}).

\subsection{Loops in $SU(2)$}\label{sec:loops}

We now compute explicitly the ring structure of $H^*_T(\Omega
SU(2); \Z)$ using the GKM graph $\Gamma\subset\mathfrak t^*$ and
the module generators $x_v$ as constructed in
Section~\ref{se:generators}. In this example, as in the previous
one, all the restrictions $x_{v}(w)$ at fixed points are
elementary tensors in $H^*_T(\{w\})\cong Sym^*(\Lambda)$. So as
before, we will represent the classes $x_v$ by drawing on every
vertex $w$  a bouquet of arrows $\beta_j\in \Lambda$ such that
$x_{v}(w)=\prod \beta_j$. The vertices with no arrows coming out
of them carry the class $0$.

\begin{figure}[h]
\begin{center}
\epsfig{figure=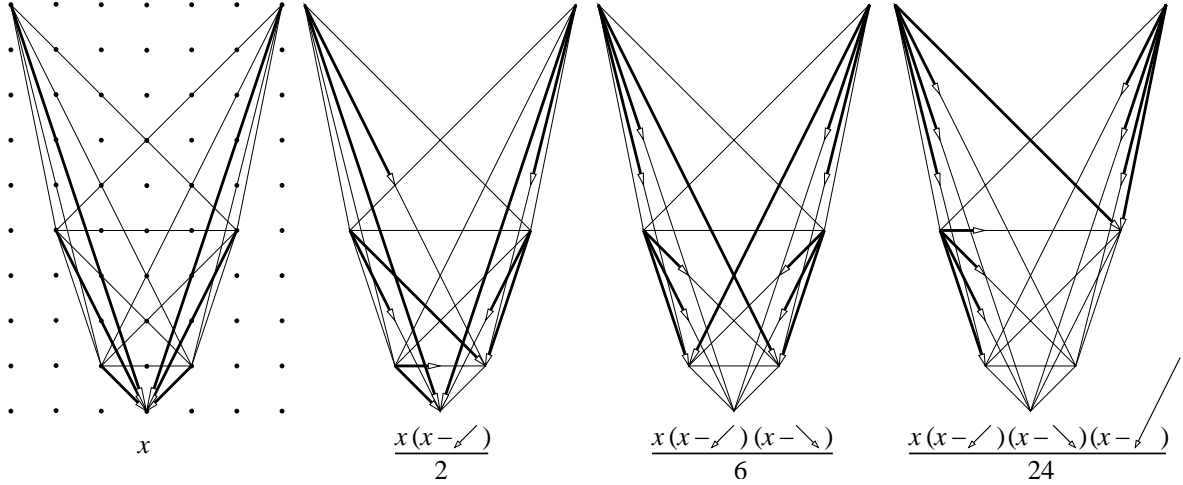,height=2.5in}
\end{center}
\begin{center}
\parbox{5in}{
\caption{ The degree $2$,$4$,$6$, and $8$ generators for
$H_T^*(\Omega SU(2); \Z)$. } \label{fig:LGfigures} }
\end{center}
\end{figure}

The first few module generators are illustrated in
Figure~\ref{fig:LGfigures}. We call $x$ the generator of degree 2,
and express the others in terms of it. The arrows in the
expressions denote elements in $H^2_T=\Lambda$.

The map $H_T^*(\Omega SU(2);\Z)\to H^*(\Omega SU(2);\Z)$ is simply
the map that sends the arrows to zero. 
So we recover the well-known fact that the ordinary
cohomology $H^*(\Omega SU(2); \Z)$ is a divided powers algebra on
a class in degree 2.

Note that the classes in Figure~\ref{fig:LGfigures} are {\em not}
generators for $K$-theory. Indeed, the conditions
(\ref{eq:iotaImage}) are only satisfied when the classes in
Figure~\ref{fig:LGfigures} are interpreted in cohomology, but not
when they are interpreted in $K$-theory.

\begin{figure}[h]
%\psfrag{1}{$1$} \psfrag{2}{$2$} \psfrag{3}{$3$}
\psfrag{x0}{$x_{1}$}\psfrag{x1}{$x_{-1}$}\psfrag{x2}{$x_2$}\psfrag{x3}{$x_{-2}$}
\begin{center}
\epsfig{figure=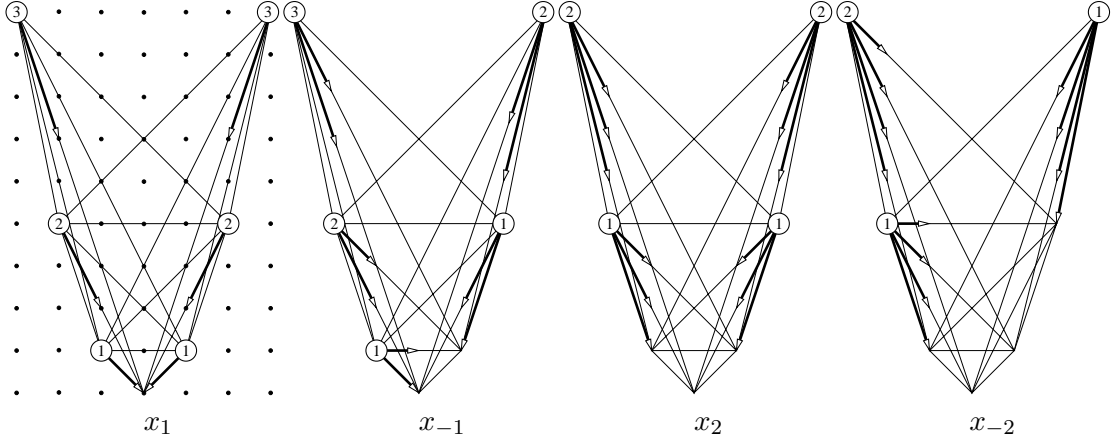,height=2.25in}
\end{center}
\begin{center}
\parbox{5in}{
\caption{ The first few module-generators of $K_T^*(\Omega
SU(2))$. (The class $x_0= 1$ is omitted.)} \label{fig:LGKthy} }
\end{center}
\end{figure}

To compute the generators of $K_T^*(\Omega SU(2))$, we introduce
the following notation. Let
\[
p_k(\lambda_1,\ldots,\lambda_n):=
(1-\lambda_1)\cdots(1-\lambda_n)\cdot
\sum_{0\le|\alpha|<k}\lambda^{\alpha},
\]
where
$\lambda^\alpha=\lambda_1^{\alpha_1}\cdots\lambda_n^{\alpha_n}$
and $|\alpha|=\alpha_1+\ldots+\alpha_n$.  The first such
polynomial $p_1(\lambda_1,\ldots,\lambda_n)$ is exactly the Euler
class $e(\bigoplus \lambda_i)$ that appeared in
Section~\ref{sec:G2}.  The other ones are slightly more
complicated. To best draw our $K$-theory classes, we introduce a
pictorial notation for $p_k(\lambda_1,\ldots,\lambda_n)$, for
$\lambda_i\in \Lambda$. We will represent them by a bouquet
$\{\lambda_i\}$ of arrows, and a small number $k$ at the vertex.
We illustrate our generators using this notation in
Figure~\ref{fig:LGKthy}. To check that these elements are indeed
the generators of $K^*_T(\Omega SU(2))$, we need to check
(\ref{xij}), which is immediate, and that they satisfy the GKM
conditions (\ref{eq:iotaImage}). These latter turn out to be quite
hard to check.

Let $a:=\put(0,3){\vector(1,0){12}}\hspace{4mm}$ and
$q:=\put(4,-3){\vector(0,1){12}}\hspace{3mm}\in K^0_T$. Let us
also identify the vertex set $F$ of $\Gamma$ with $\Z$ by taking
the horizontal coordinate.  The class $x_i$ drawn in
Figure~\ref{fig:LGKthy} is given by
\[
x_i(m)=
\begin{cases}
p_{m-k}\big(a^{-1}q^{-m-k},a^{-1}q^{-m-k+1},\ldots,a^{-1}q^{-m+\ell}\big)
& \mbox{if } m>k \\
0 & \mbox{if } -\ell \leq m\leq k\\
p_{-m-\ell}(aq^{m-\ell},aq^{m-\ell+1},\ldots,aq^{m+k}) & \mbox{if
} m<-\ell,
\end{cases}
\]
where  $\ell = |i|-1$ and $k = | i - \frac{1}{2}| - \frac{1}{2}$.
Given an edge $(m,n)\in\Gamma$, we must check the condition given
in (\ref{eq:iotaImage}), namely that the Euler class $1-aq^{m+n}$
divides the difference
$$
x_i(m)-x_i(n).
$$
This involves several different cases. However, the problem has a
few symmetries that allow us to reduce the cases to the following
three.

If $m$ is between $-\ell$ and $k$ then $x_i(n)$ has either
$(1-aq^{m+n})$ or $(1-a^{-1}q^{-m-n})$ as a factor and we are
done.

If both $m$ and $n$ are bigger than $k$, then we must check that
$1-aq^{m+n}$ divides
\begin{equation}\label{eq:divisi}
p_{m-k}\big(a^{-1}q^{-m-k},\ldots,a^{-1}q^{-m+\ell}\big)-p_{n-k}\big(a^{-1}q^{-n-k},\ldots,a^{-1}q^{-n+\ell}\big).
\end{equation}
This is equivalent to checking that \eqref{eq:divisi} evaluates to
0 after setting $a^{-1}=q^{m+n}$. So we are reduced to checking
that
\[
p_{m-k}(q^{n-k},\ldots,q^{n+\ell})=p_{n-k}(q^{m-k},\ldots,q^{m+\ell}).
\]
The above formula is invariant under adding the same constant to
the indices $m,$ $n$ and $k$, and subtracting it from $\ell$. So
by letting $k=0$, we must prove the equivalent formula
\begin{equation}\label{m=n}
p_{m}(q^{n},\ldots,q^{n+\ell})=p_{n}(q^{m},\ldots,q^{m+\ell}).
\end{equation}
This is the content of Lemma \ref{lemma:qqq}.

Finally, if $m>k$ and $n<-\ell$ then we are reduced to checking
that
$$
p_{m-k}(q^{n-k},\ldots,q^{n+\ell})=p_{-n-\ell}(q^{-m-\ell},\ldots,q^{-m+k}).
$$
By replacing $q$ with $q^{-1}$, reversing the order of the
arguments in the polynomial $p$, and a couple changes of indices,
this also reduces to Lemma \ref{lemma:qqq}.

\begin{lemma}\label{lemma:qqq}
The expression
$$
a_{mn\ell}:=p_m(q^n,q^{n+1},\ldots,q^{n+\ell})
$$
is symmetric in $m$ and $n$.
\end{lemma}

\begin{proof}
Let $\binom{}{}_q$ denote the quantum binomial coefficient
$$
{a\choose b}_{\!q}=\frac{a!_q}{b!_q(a-b)!_q},
$$
where $a!_q$ is the $q$-factorial\footnote{Some authors define the
quantum factorial $a!_q$ to be $1(1+q)(1+q+q^2)\cdots (1+q+\cdots
+q^a)$. This agrees with our expression up to a power of $1-q$.}
$a!_q=(1-q)(1-q^2)\ldots(1-q^a)$. We can then rewrite the
expression $a_{mn\ell}$ as
\begin{equation}\label{eq:qbinom}
a_{mn\ell}=(1-q^n)\cdots(1-q^{n+\ell})\cdot\sum_{i=0}^{m-1}q^{in}{\ell+i\choose
\ell}_{\!q}.
\end{equation}
See for example \cite[\S 3.3]{Andrews} for more detail. In
particular, \eqref{eq:qbinom} is a truncated version of
Equation~(3.3.7) in \cite{Andrews}.

Now recall from \cite{Zei93} that a ``difference form'' $$\omega =
f(i,j)\delta i + g(i,j)\delta j$$ has ``exterior difference''
$$
d\omega = \big[f(i,j+1)-f(i,j)\big]\delta j\ \delta i + \big[
g(i+1,j)-g(i,j)\big] \delta i\ \delta j,
$$
where $\delta i$ and $\delta j$ are anti-commuting symbols. Such a
difference form can be viewed as a cellular $1$-cochain on the
standard square tiling of $\R^2$,  the exterior difference being
the usual cellular coboundary operator.  Consider the difference
form
$$
\omega = q^{ij}\frac{(i+\ell)!_q(j+\ell)!_q}{i!_qj!_q\ell!_q}
\big[ (1-q^j)\delta i + (1-q^i)\delta j\big].
$$
It is an easy exercise to verify that $\omega$ is closed.
Therefore, by the discrete Stokes' theorem \cite{Zei93},
$$
\int_{\partial L} \omega = 0,
$$
where $L$ is the rectangle  $[0,m]\times [0,n]$.  One now checks
that the above integral is zero on the sides $\{ 0\} \times [0,n]$
and $[0,m]\times\{ 0\}$, and equals $a_{nm\ell}$ and $-a_{mn\ell}$
on the remaining two sides.
\end{proof}

\begin{remark}
We do not know whether the generators illustrated in
Figure~\ref{fig:LGKthy} are the same as those mentioned in
Remark~\ref{rem:KthGen}.
\end{remark}

\subsection{A homogeneous space of type $A_1^{(4)}$}\label{twLoops}

For our last example, we let $\G$ be the affine group associated
to the Cartan matrix
\[
\left[
\begin{array}{rr} 2 & -1 \\ -4 & 2 \end{array} \right].
\]
This group is \(\widehat{LSL(3,\C)}^{\Z/2\Z} \rtimes \C^*,\) where
the $\Z/2\Z$-action on $LSL(3,\C)$ is given by precomposition with
the antipodal map \(z \mapsto -z\) on $\C^*$ and composition with
the outer automorphism \(A \mapsto (A^t)^{-1}\) of $SL(3,\C)$.

\begin{figure}[h]
\begin{center}
\epsfig{figure=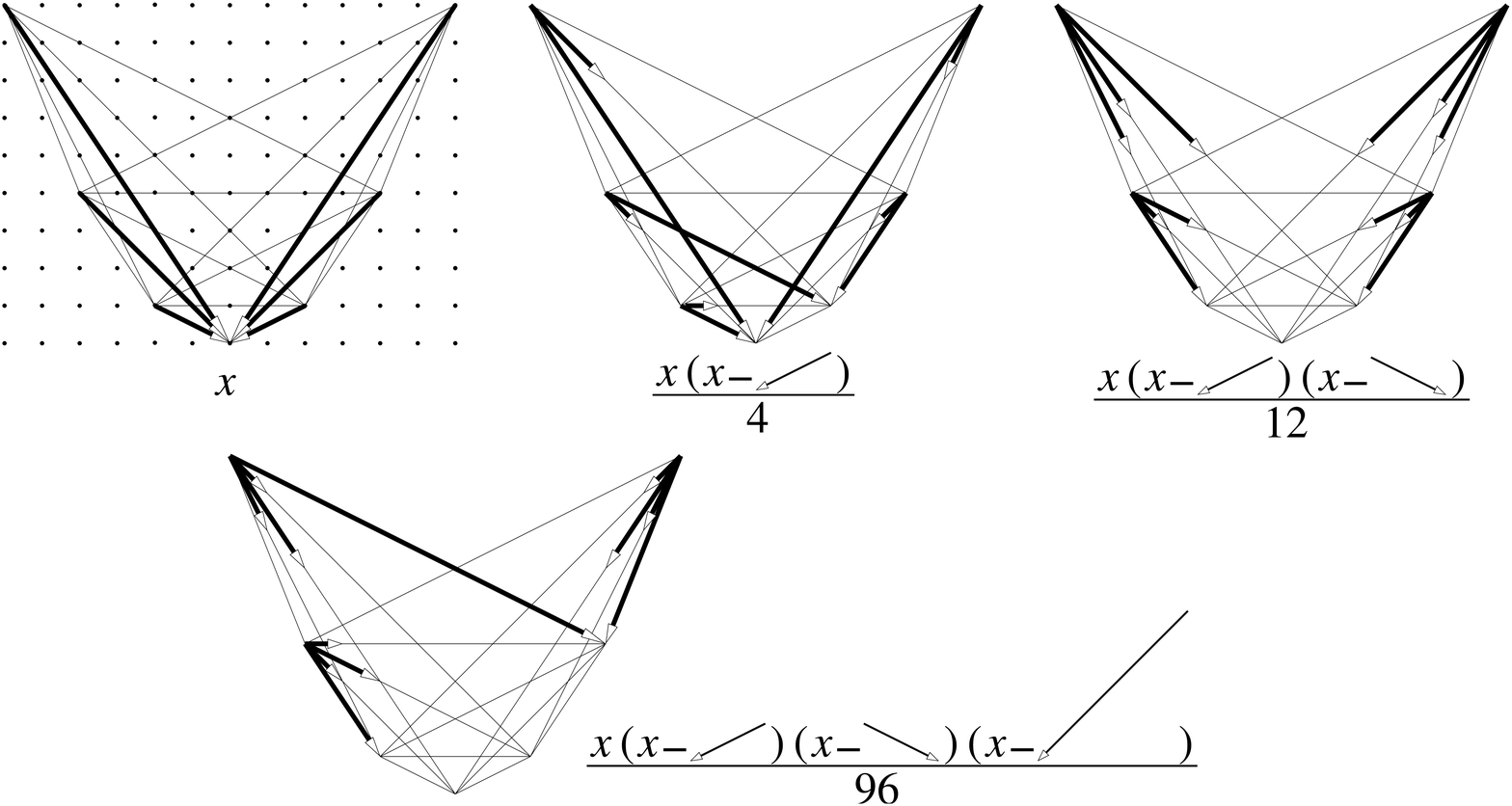,height=3.25in}
\end{center}
\begin{center}
\parbox{5in}{
\caption{ The degree 2, 4, 6, and 8 generators for $H^*_T(\G/\P;
\Z)$.}\label{fig:twLGfigures} }
\end{center}
\end{figure}

We consider the homogeneous space $\G/\P$ where the parabolic $\P$
has Lie algebra generated by $\mathfrak{b}$ and the negative of
the simple short root. The degree 2, 4, 6, and 8 module generators
for $H_T^*(\G/\P;\Z)$ are illustrated in
Figure~\ref{fig:twLGfigures}. The denominator in the degree $n$-th
module generator is given by \(n!2^{\lfloor n/2 \rfloor}.\)

\bibliographystyle{alpha}
\bibliography{ref}

%\end{spacing}

\end{document}